\documentclass[12pt]{article}
\usepackage{latexsym, amsfonts, amsmath, amssymb}
\usepackage{graphicx}
\usepackage{graphicx}
\usepackage{amsthm}
\usepackage{yfonts}
\usepackage{epsfig}
\usepackage[utf8]{inputenc}
\usepackage[english]{babel}
\usepackage{empheq}

\usepackage{fullpage}

\usepackage{color}
\definecolor{purple}{rgb}{.8,0,.8}

\def\boh#1{{\bf #1} \ }

\begin{document}

\def\baselinestretch{1.2}

\def\ostrut#1#2{\hbox{\vrule height #1pt depth #2pt width 0pt}}

\setlength{\parskip}{12pt}

\begin{titlepage}

\vglue 1in
\begin{center}
{\Large \bf Points of constancy of the periodic\ostrut0{11} linearized Korteweg--deVries equation}
\vskip 0.5cm

Peter J.~Olver$^{1,a}$ and Efstratios Tsatis$^{2,b}$

\vspace{0.1cm}

{\it ${}^1$School of Mathematics, University of Minnesota, Minneapolis, MN \ 55455, USA}

\vspace{-0.1cm}

{\it ${}^2$ 8 Kotylaiou Street, Athens 11364, Greece }

\vspace{0.1cm}

{\tt  ${}^a$olver@umn.edu,  ${}^b$efstratiostsatis@hotmail.com} \\

\medskip

\end{center}
\vskip .5in
\centerline{\bf Abstract}

\vskip20pt
We investigate the points of constancy in the piecewise constant solution profiles  of the periodic linearized Korteweg--deVries equation with step function initial data at rational times.  The solution formulas are given by certain Weyl sums, and we employ number theoretic techniques, including Kummer sums, in our analysis. These results constitute an initial  attempt to understand the complementary phenomenon of ``fractalization'' at irrational times.

\vspace{0.2cm}

\end{titlepage}


\renewcommand\Re{\operatorname{Re}}
\renewcommand\Im{\operatorname{Im}}

\newtheorem{remark}{Remark}
\newtheorem{prop}{Proposition}

\let\rf\cite
\def\blackbox{\vrule height7pt width5pt depth0pt}

\def\Eq#1$$#2$${\StEq#1  \EnEq{#2}}
\def\StEq#1 #2\EnEq#3{\begin{equation}\label{#1eq*} #3\end{equation}}

\def\eqrefz#1{\ref{#1eq*}}
\def\eq#1{{\rm (\eqrefz{#1})}}
\def\eqe{equation~\eq} \def\eqse{equations~\eq} \def\eqE{Equation~\eq} 
\def\eqf{formula~\eq} \def\eqF{Formula~\eq} \def\eqfe{formulae~\eq}
\def\eqs#1#2{\eq{#1},~\eq{#2}}
\def\eqas#1#2{\eq{#1} and~\eq{#2}}

\newdimen\eqjot~\eqjot = 1\jot
\def\openupeq{\openup \the\eqjot}
\newdimen\tablejot \tablejot = 1\jot
\def\openuptable{\openup \the\tablejot}

\def\crh#1{{}\cr\hskip#1pt{}}
\def\crhz{\crh0} \def\crhx{\crh{10}} \def\crhq{\crh{20}} \def\crhQ{\crh{40}}
\def\crr#1#2{#2{}\cr\hskip#1pt{}#2}
\def\crz{\crr0} \def\crx{\crr{10}} \def\crq{\crr{20}} \def\crQ{\crr{40}}
\def\cth#1#2{#2{}\cr&\hskip#1pt{}#2}
\def\ctz{\cth0} \def\ctx{\cth{10}} \def\ctq{\cth{20}} \def\ctQ{\cth{40}}
\def\creq{\ctq}\def\crpm{\crq}
\def\cthn#1#2{{}\cr&\hskip#1pt{}#2}
\def\ctnz{\cthn0} \def\ctnx{\cthn{10}} \def\ctnq{\cthn{20}} \def\ctnQ{\cthn{40}}
\def\creqn{\ctnq}
\def\crs{\cr\noalign{\medskip}}

\def\addtab#1={#1\;&=}
\def\addtabe#1=#2={#1=#2\;&=}

%
%
%
%
%

\def\aaeq#1#2{{\def\\{\cr}\vcenter{\openupeq \halign{$\displaystyle 
   \hfil##\;$&$\displaystyle##\hfil$&&\hskip#1pt$\displaystyle##\hfil$\cr#2\cr}}}}
\def\aeq{\aaeq{20}}

\def\qaeq#1#2{{\def\\{&}\vcenter{\openupeq\halign{$\displaystyle
   ##\hfil$&&\hskip#1pt$\displaystyle##\hfil$\cr #2\cr}}}}
\def\qqeq{\qaeq{40}} \def\qeqq{\qaeq{40}}
\def\req{\qaeq{30}} \def\qxeq{\qaeq{30}} \def\qeqx{\qaeq{30}}
\def\qeq{\qaeq{20}} \def\weq{\qaeq{15}}
\def\xeq{\qaeq{10}} \def\yeq{\qaeq{7}} \def\zeq{\qaeq{3}}

\def\caeq#1#2{{\def\\{\cr}\vcenter{\openupeq \halign
    {$\hfil\displaystyle ##\hfil$&&$\hfil\hskip#1pt\displaystyle ##\hfil$\cr#2\cr}}}}
\def\ceq{\caeq{20}} \def\deq{\caeq{10}} \def\doeq{\caeq{3}}

\def\ezeq#1#2#3{{\def\\{\cr#1}\vcenter{\openupeq \halign{$\displaystyle 
   \hfil##$&$\displaystyle##\hfil$&&\hskip#2pt$\displaystyle##\hfil$\cr#1#3\cr}}}}
\def\eaeq{\ezeq\addtab} \def\eeeq{\ezeq\addtabe{20}}
\def\eqeq{\eaeq{40}} \def\eeqq{\eaeq{40}} 
\def\exeq{\eaeq{30}} \def\eeqx{\eaeq{30}}
\def\eeq{\eaeq{20}} \def\feq{\eaeq{10}} \def\heq{\eaeq{3}}

\def\iaeq#1#2#3{{\def\\{\cr\hskip#1pt}\vcenter{\openupeq\halign{$\displaystyle
   ##\hfil$&&\hskip#2pt$\displaystyle##\hfil$\cr #3\cr}}}}
\def\ibeq#1{\iaeq{#1}{#1}}
\def\iqeq{\iaeq{40}{40}}  \def\iQeq{\iaeq{80}{80}}
\def\ieq{\iaeq{20}{20}} \def\jeq{\iaeq{10}{10}} \def\keq{\iaeq{10}{3}}

\def\saeq#1#2{{\def\\{\cr}\vcenter{\openup1\jot \halign{$\displaystyle
   ##\hfil$&&\hskip#1pt$\displaystyle##\hfil$\cr #2\cr}}}}
\def\sqeq{\saeq{40}} \def\seqq{\saeq{40}}
\def\sxeq{\saeq{30}} \def\seqx{\saeq{30}}
\def\seq{\saeq{20}} \def\teq{\saeq{10}} \def\ueq{\saeq{3}}

\def\macases#1#2{\left\{\enspace\saeq{#1}{#2}\right.}
\def\mcases{\macases{20}}
\def\xcases{\macases{15}}
\def\tcases{\macases{10}}
\def\ucases{\macases{3}}

\def\parens#1{({#1})}                       \let\pa\parens
\def\parenz#1{\bigl({#1}\bigr)}             \let\paz\parenz
\def\bparens#1{\bigl(\,{#1}\,\bigr)}        \let\bpa\bparens
\def\Parens#1{\left (\,{#1}\,\right )}      \let\Pa\Parens
\def\Parenz#1{\left ({#1}\right )}          \let\Paz\Parenz
\def\brackets#1{[\mskip2mu{#1}\mskip 2mu]}  \let\bk\brackets 
\def\bracketz#1{\bigl [{#1}\bigr ]}         \let\brz\bracketz \let\bkz\bracketz
\def\bbrackets#1{\bigl [\,{#1}\,\bigr ]}    \let\bbr\bbrackets \let\bbk\bbrackets
\def\Brackets#1{\left [\,{#1}\,\right ]}    \let\Br\Brackets \let\Bk\Brackets
\def\Bracketz#1{\left [{#1}\right ]}        \let\Brz\Bracketz \let\Bkz\Bracketz 
\def\braces#1{\{\mskip2mu{#1}\mskip 2mu\}}  \let\bc\braces
\def\bracez#1{\bigl \{{#1}\bigr \}}         \let\bcz\bracez
\def\bbraces#1{\bigl \{\,{#1}\,\bigr \}}    \let\bbc\bbraces
\def\Braces#1{\left \{\,{#1}\,\right \}}    \let\Bc\Braces
\def\Bracez#1{\left \{{#1}\right \}}        \let\Bcz\Bracez

\newtheorem{Theorem}{Theorem}
\newtheorem{Corollary}[Theorem]{Corollary}
\newtheorem{Conjecture}[Theorem]{Conjecture}
\newtheorem{Definition}[Theorem]{Definition}
\newtheorem{Lemma}[Theorem]{Lemma}
\newtheorem{Example}[Theorem]{Example}
\newtheorem{Exercise}[Theorem]{Exercise}
\newtheorem{Proposition}[Theorem]{Proposition}
\newtheorem{Remark}[Theorem]{Remark} 
\newtheorem{Algorithm}{Algorithm}[section]

\def\Ex#1{\exdeclare{Example}{#1}} \def\ex{\stwrite{Example}}
\def\exs#1#2{Examples \stn{#1} and~\stn{#2}}
\def\Ez#1{\exdeclare{Exercise}{#1}} \def\ez{\stwrite{Exercise}}
\def\ezs#1#2{Exercises \stn{#1} and~\stn{#2}}
\def\Th#1{\stdeclare{Theorem}{#1}} \def\th{\stwrite{Theorem}}
\def\ths#1#2{Theorems \stn{#1} and~\stn{#2}}
\def\Lm#1{\stdeclare{Lemma}{#1}} \def\lm{\stwrite{Lemma}}
\def\lms#1#2{Lemmas \stn{#1} and~\stn{#2}}
\def\Pr#1{\stdeclare{Proposition}{#1}} \def\pr{\stwrite{Proposition}}
\def\prs#1#2{Propositions \stn{#1} and~\stn{#2}}
\def\Co#1{\stdeclare{Corollary}{#1}} \def\co{\stwrite{Corollary}}
\def\cors#1#2{Corollarys \stn{#1} and~\stn{#2}}
\def\Df#1{\exdeclare{Definition}{#1}} \def\df{\stwrite{Definition}}
\def\dfs#1#2{Definitions \stn{#1} and~\stn{#2}}
\def\Cj#1{\stdeclare{Conjecture}{#1}} \def\cj{\stwrite{Conjecture}}
\def\cjs#1#2{Conjectures \stn{#1} and~\stn{#2}}

\def\exdeclare#1#2#3\par{\begin{#1}\label{#2}{\rm #3}\end{#1}\rm}
\def\stdeclare#1#2#3\par{\begin{#1}\label{#2}{#3}\end{#1}}

\def\stwrite#1#2{#1 \stn{#2}}  
\def\stn{\ref}

\def\Rmk#1\par{\is{Remark}: \ #1\par\medbreak}

\newif\ifnames \namesfalse
\def\namez#1;{\namesfalse\namezz#1,;}
\def\namezz#1,#2,#3;{\ifx,#3,\ifnames \ and\fi\fi #1,#2,\ifx,#3,\else
    \namestrue\namezz #3;\fi}

\def\key#1 #2\par{\bibitem{#1} #2\par}

\def\book#1;#2;#3\par{\namez#1;{\it #2},#3.}
\def\paper#1;#2;#3; #4(#5)#6\par{\namez#1;#2,{\it #3}\if a#4, to appear.\else
     { \bf #4}(#5),#6.\fi}
\def\inbook#1;#2;#3;#4\par{\namez#1;#2, {\sl in\/}:{\it#3},#4.}
\def\appx#1;#2;#3;#4\par{\namez#1;#2, {\sl appendix in\/}:{\it #3},#4.\par }
\def\thesis#1;#2;#3\par{\namez#1;{\it #2}, Ph.D.~Thesis,#3.}
\def\uthesis#1;#2;#3\par{\namez#1;{\it #2}, Ph.D.~Thesis, University of Minnesota,#3.}
\def\preprint#1;#2;#3\par{\namez#1;#2, preprint,#3.}
\def\upreprint#1;#2;#3\par{\namez#1;#2, preprint, University of Minnesota,#3.}
\def\prepare#1;#2;#3\par{\namez#1;#2, in preparation.}
\def\personal#1;#2\par{\namez#1; personal communication,#2.}
\def\submit#1;#2;#3\par{\namez#1;#2, submitted.}
\def\other#1\par{#1.}

\def\ro#1{{\rm #1}} \def\rbox#1{\hbox{\rm #1}}
\def\roh{\rbox}
\def\roz#1{\  \rbox{#1}\ }
\def\rox#1{\quad \rbox{#1}\quad }
\def\roy#1{\ \ \rbox{#1}\ \ }
\def\roq#1{\qquad \rbox{#1}\qquad }
\def\roqx#1{\qquad \quad \rbox{#1}\quad \qquad }
\def\roqq#1{\qquad \qquad \rbox{#1}\qquad \qquad }

\def\slbox#1{\hbox{\it #1}}
\def\soq#1{\qquad \slbox{#1}\qquad }

\def\where{\roq{where}}

\def\Bbb#1{{\mathbb#1}}

\def\R{\Bbb R}  \def\C{\Bbb C}
\def\N{\Bbb N}  \def\P{\Bbb P}
\def\Q{\Bbb Q}  \def\Z{\Bbb Z}

\def\hexnumber#1{\ifcase#1 0\or1\or2\or3\or4\or5\or6\or7\or8\or9\or
 A\or B\or C\or D\or E\or F\fi}
\edef\msbhx{\hexnumber\symAMSb}

\mathchardef\emptyset="0\msbhx3F
\mathchardef\subsetneq="2\msbhx28
\mathchardef\supsetneq="2\msbhx29

\def\chat{\widehat c\:}
\def\qhat{\widehat q\:}

\def\msk#1{\mskip #1 mu} \def\:{\mskip2mu}
\def\pii{\: \pi} \def\i{\,\ro i\,}

\def\rg#1#2{#1=1,\ldots,#2} \def\rgo#1#2{#1=0,\ldots,#2} 

\def\PC{\mathfrak{P}}
\def\q{\mathfrak{q}} \def\qi{\q_i} 
\def\qh{\widehat \q} \def\qhi{\qh_i}
\def\elt{\widetilde \ell}

\def\Sumu#1{\sum _{#1}\>}
\def\Summ#1#2#3{\sum _{#1\:=\:#2}^{#3}\>}
\def\Sum#1#2{\Summ{#1}1{#2}} \def\Sumo#1#2{\Summ{#1}0{#2}} 

\def\operator#1{\expandafter\def\csname#1\endcsname{\mathop{\rm #1}\nolimits}}

\def\fracz#1#2{#1/#2}
\def\dsty{\displaystyle} 

\def\Se{S_e} \def\So{S_o}

\def\modx{{\rm mod} \>}

\numberwithin{equation}{section}

\section{Introduction}

The starting point of our investigations is the initial-boundary value problem for the most basic linearly dispersive wave equation for a function $u(t,x)$,
\Eq{lkdv}
$$\frac{\partial u}{\partial t}=\frac{\partial^3 u}{\partial x^3},$$
known as the \emph{linearized Korteweg--deVries} or {\emph{Airy partial differential equation} because its fundamental solution can be expressed in terms of the Airy function, \rf P}.  We are interested in the case of periodic boundary conditions 
\Eq{pbc}
$$u(t,0)=u(t,2\pii), \qquad  \frac{\partial u}{\partial x}(t,0)=\frac{\partial u}{\partial x}(t,2\pii), \qquad  \frac{\partial^2 u}{\partial x^2}(t,0)=\frac{\partial^2 u}{\partial x^2}(t,2\pii) 
$$
on the interval $0 \leq x \leq 2\pii$.  The initial data of interest is provided by the discontinuous unit step function:
\Eq{step}
$$u(0,x)=\sigma(x)\equiv \left\{
     \begin{array}{ll}
       0,  & \quad 0<x<\pi,  \\
         1, & \quad \pi<x<2\pii, \\
     \end{array}  
   \right. 
$$
{sometimes referred to as the \emph{Riemann problem}, which is of fundamental importance in the study of hyperbolic wave equations and shock waves, \rf{Smoller}}. The precise value assigned at its discontinuity is not important, although choosing $\sigma (x) = \frac{1}{2}$ at $x=0,\pi,2\pii$ is consistent with Fourier analysis, \rf P.  
The boundary conditions allow one to 
extend the initial data and solution to be $2\pii$-periodic functions in $x$, with jumps of magnitude $\pm\:1$ at integer multiples of $\pi$, and we use the same notation, so that $\sigma (x)$ denotes the $2\pii$-periodic extension of the unit step function throughout.

It was shown in~\rf{Odq} {that, in striking contrast to the smooth evolution of the Riemann solution on the line,} the resulting (weak) solution to the periodic initial-boundary value problem~\eq{lkdv}--\eq{step} exhibits \emph{dispersive quantization}, also known as the \emph{Talbot effect}, \rf{BerryKlein,BMS,ErTz}.  Namely, whenever 
time is a rational multiple of $\pi$, the solution profile is discontinuous, but piecewise constant. At irrational times the solution profile is continuous, but fractal and non-differentiable, \rf{Osk1,ChErTz,ErSh}.  In this paper, we 
concentrate on piecewise constant profiles at rational times, since these are also of interest in number theory, being closely related to Weyl sums~\rf{Vino}. This surprising phenomenon is observed in other 
linearly dispersive partial differential equations, such as the linear Schr\"odinger equation, and also in nonlinear equations, both integrable and non-integrable, including the nonlinear Schr\"odinger and Korteweg--deVries equations with a variety of nonlinearities, \rf{COdisp}.  Experimental confirmations of the Talbot effect in both optics and atoms are described in~\rf{BMS}.
 
The precise theorem concerning the above initial-boundary value problem for the Airy equation~\eq{lkdv} can be stated as follows:

\Th1
Let $p/q \in \Q$ be a rational number with $p$ and $q$ having no common factors.  
Then the solution  to the initial-boundary value problem  at time $t = \pi\msk1 p/q$ is constant on every subinterval $\pi j/q < x < \pi(j+1)/q$ for $j =  0, \ldots,2\:q-1$.

Thus, at rational time $t = \pi\msk1 p/q$ (relative to the length $2\pii$ of the interval), the solution achieves a constant value on each specified subinterval of length $\pi/q$.  It was further noted 
that, often, the solution is in fact constant on longer subintervals than those specified in \th1, and the question arises of how to characterize these ``regions of constancy".  This problem and its number-
theoretic implications form the focus of this paper.

It was further proved in \rf{Odq} that the \emph{fundamental solution} $u = F(t,x)$ to the periodic initial-boundary value problem, meaning the one with initial conditions a (periodically extended) delta 
function, $F(0,x) = \delta (x)$, is, at the rational time $t = \pi\msk1 p/q$, a finite linear combination of (periodically extended) delta functions, often called a \emph{Dirac comb}, based at the rational 
nodes  $\pi\msk1\ell/q$, so
\Eq{fs}
$$F(\pi \msk1 p/q,x) = \Sumo l{2q-1} \beta _\ell(p/q) \,\delta (x - \pi\msk1\ell/q).$$
{The coefficients $\beta _\ell(p/q)$ turn out to be expressible in terms of \emph{Kummer sums}, \eq{TS}, whose formulae are given in \eq{betal} below.  The fact that, at rational times, the fundamental solution can be represented by  such sums underlies the deep connections between number theory and the solutions to elementary linear partial differential equations on periodic domains that will be explored in this paper.} 
This fact also underlies the phenomenon of \emph{quantum revival},\rf{BerryKlein}, in which an initially concentrated wave function, representing, say, an electron in an atomic orbit, at first spreads out but subsequently relocalizes at all rational times, the number of localization sites depending upon the size of the denominator $q$.   

{A remarkable consequence of formula~\eq{fs} is} that the value \emph{any} $2\pii$ periodic (in $x$) solution at each rational time is a \emph{linear combination of finitely many translates of its initial profile} $u(0,x) = f(x)$:
\Eq{utrans}
$$u(\pi \msk1 p/q,x) = \Sumo \ell{2q-1} \beta _\ell(p/q) \,f(x - \pi\msk2 \ell/q).$$
Thus, at rational times, the value of the solution at a point $x$ \emph{depends on only finitely many values of its initial data}!
In particular, for the periodically extended step function initial data~\eq{step}, the superposition formula implies that the resulting solution at rational time is given by
\Eq{steptrans}
$$u(\pi\msk1  p/q,x) = \Sumo \ell{2q-1} \beta _\ell(p/q) \,\sigma (x - \pi\msk2 \ell/q).$$

Since the delta function is even, $\delta (-\:x) = \delta (x)$, the same holds for the fundamental solution, and hence
\Eq{evenc}
$$\beta _{2q-\ell}(p/q)= \beta _\ell(p/q).$$ 
Furthermore, substituting the constant solution $u(t,x)~\equiv 1$ into~\eqf{utrans}, we deduce that
\Eq{csum}
$$\Sumo \ell{2q-1} \beta _\ell(p/q)  = 1.$$

For any piecewise continuous function, let
$$\Delta f(x) = f(x^+) - f(x^-)$$
denote its \emph{jump} at the point $x$; thus $\Delta f(x) = 0$ if $f$ is continuous at $x$.
In particular, 
$$\req{\Delta \sigma (x) = \mcases{1,& x = (2\:j+1) \pii,\\-1& x = 2\:j\pii,\\0,&\roh{otherwise,}}\\j \in \Z.}$$
{In view of \eq{evenc}}, the jumps of the solution~\eq{steptrans} at the nodes are
\Eq{jumpu}
$$\xeq{\Delta u(\pi \msk1 p/q,\pi \msk1 j/q) = \beta _{j+q}(p/q) - \beta _j(p/q) = -\,\Delta u(\pi \msk1 p/q,\pi \msk1 (j+q)/q),\\ \rgo j{q-1}.}$$
In other words, the jumps are symmetrically paired with jumps of opposite magnitude around the midpoint $x = \pi$ or, equivalently, around $j=q$.

More specifically, we call a connected open subinterval $I\subset \R$ a \emph{region of constancy} for the specified rational time $t = \pi \msk1p/q$ if $u(t,x)$ is constant for all $x \in I$.  As an immediate 
consequence of \th1, we deduce that every region of constancy has the form
\Eq{rc}
$$ \pi \msk1 j/q < x < \pi\msk1  k/q \roq{for} j,k \in \Z \rox{with} 1 \leq k-j < 2 \:q.$$
Note that the strict inequality on $k-j$ is because the solution is clearly not constant on the entire interval $[0,2\pii]$ as otherwise it would be the trivial constant solution.  Let us call an integer 
$m$ a \emph{point of constancy} if $x_m = \pi \msk1 m/q$ belongs to a region of constancy.  Thus, the region of constancy given in~\eq{rc} contains the $k-j-1$ points of constancy ${j+1},{j+2},\ldots,{k-1}$. On 
the other hand, when $k-j=1$, which is the minimum length of a region of constancy, {the region} contains no points of constancy.  Given a rational time $t = \pi \msk1 p/q$, the aim of this paper is to characterize all 
points of constancy, which is equivalent to characterizing all regions of constancy of the solution $u(t,x)$. As we will see, and as anticipated in~\rf{Odq}, this characterization will be number theoretic in 
nature, and in fact relies on the properties of (partial) Kummer sums~\rf{HBPa}.

More explicitly, let
\Eq{aj}
$$a_j = a_j\left(\frac{p}{q}\right) = u\left(\frac{\pi \msk1  p}{q},x\right) \roq{for any} \frac{\pi \msk1  j}q < x < \frac{\pi(j+1)}q$$ 
denote the constant value of the solution on the indicated subinterval at time $t = \pi\msk1  p/q$.  Clearly, $j$ is a point of constancy if and only if 
\Eq{pc}
$$\boxed{a_j = a_{j-1}.}$$
This is the key equation we will analyze in what follows.  Referring to~\eq{steptrans}, we see that
\Eq{ajcj}
$$a_j\left(\frac{p}{q}\right) = \mcases{\Summ i{j+1}{j+q} \beta _j(p/q),& 0 \leq j < q,\\ \Sumo i{j-q} \beta _j(p/q) \; + \; \Summ i{j+1}{2q-1} \beta _j(p/q),& q \leq j < 2\:q.}$$

\Df{PC}
The \emph{set of all points of constancy} at $t=\pi \msk1  p/q$ is denoted  
\Eq{PC}
$$\PC\left(t=\frac{\pi \msk1 p}{q}\right)\equiv\PC\left(\frac{p}{q}\right)=\left\{j\in[0,2\:q-1]\cap\Z\ \bigg|\  a_{j}-a_{j-1}=0\right\}.$$

Our goal is to explicitly characterize the set~\eq{PC}.  A key result is that the set $\PC(p/q)$ will be given in terms of the primes $\qi$ appearing in the prime factorization of $q$, which we denote by
\Eq{qfac}
$$q=\prod_{i=1}^{m}\>\qi^{n_i}, $$
where the prime factors $\qi$ are distinct and the powers $n_i \geq 1$.
We find that the set $\PC(p/q)$ equals the union of the corresponding sets of points of constancy at the times given by the prime powers appearing in the factorization, so
\Eq{PCfac}
$$\PC\left(\frac{p}{q}\right)=\bigcup_{i} \; \PC\left(\frac{p}{\qi^{n_i}}\right).$$
This observation allows us to study each prime separately. It is important to note that the resulting sets depend on both the prime factor $\qi$ and the numerator $p$. 

This paper is organized as follows.  After some preliminary computations, in Section~\ref{sec:squarefree} we introduce the notion of permutation polynomials in order to analyze the case when $q$ is square free, i.e., all $n_i = 1$ in~\eq{qfac}. In Section~\ref{sec:qodd} we discuss the case when $q$ is odd, and in Section~\ref{sec:qeven} the case when $q$ is even. {Our main results concerning the points of constancy can be found in the general structural decompositions \eq{49}, \eq{49e}, and the various cases \eq{oddqc}, \eq{oddq2}, \eq{oddq3}, \eq{even2}, \eq{even4}, \eq{even8}, \eq{evenqn}, \eq{evenq2}.  This is not an exhaustive list, and while we have some further experimental data for some of the more complicated cases not covered, a complete resolution of this problem will necessitate a more detailed analysis, perhaps requiring more powerful number-theoretic tools. Finally,} in Section~\ref{sec:conclusion}, we give our concluding remarks and outline some directions for  future investigations.

\section{Preliminaries}\label{sec:prelim}

{The complex Fourier series of a $2\pii$ periodic function $f(x)$ is written as
\Eq{Fs}
$$f(x)\sim\sum_{k-\infty}^{\infty}c_ke^{\i k\:x},$$
where $c_k$ are its Fourier coefficients, \rf P, and, in conformity with the conventions of Fourier analysis, we use~$\sim$ rather than~$=$ to indicate that the Fourier series is formal and, without additional assumptions or analysis, its convergence is not  guaranteed.  A straightforward computation expresses the solution to the periodic initial-boundary value problem~\eq{lkdv} for the Airy equation as a time-dependent Fourier series:
\Eq{Asol}
$$u(t,x)\sim\sum_{k-\infty}^{\infty}b_ke^{\i (kx-k^3t)},$$
where, the $b_k$ are the Fourier coefficients of the initial data $u(0,x)$, which in the case of the step function~\eq{step}, are}
\begin{equation}
   b_{k} = \mcases{\frac{\i}{\pi k},  & k \roh{ odd},  \\
      \frac{1}{2},  & k=0,     \\
               0, & 0\neq k \roh{ even}. } 
\end{equation} 
Inserting these particular values into the solution formula \eq{Asol}, and rewriting the result in terms of real trigonometric functions, we obtain the (formal) Fourier expansion of the solution to the 
original initial-boundary value problem
\begin{equation}
u(t,x)\;\sim \;\frac{1}{2}\;-\;\frac{2}{\pi}\;\sum_{j=0}^{\infty}\frac{\sin\bbk{(2j+1)x-(2j+1)^3 t}}{2j+1} .
\end{equation}
{The resulting Fourier series can be shown to be conditionally convergent, and represents a weak solution to the original initial-boundary value problem.}

\th1 assures us that at a rational time, the solution assumes the form
\Eq{upwc}
$$u\Pa{\frac{\pi\msk1 p}q,x}=\sum_{j=0}^{2q-1}\> a_{j}\sigma^{j,q}(x),$$
where, as in~\eq{aj}, $a_j $ denotes the constant value of $u(\pi  p/q,x)$ on the subinterval $\pi  j /q < x < \pi(j+1)/q$, while
\begin{equation}
\sigma^{j,q}(x)=\left\{
     \begin{array}{ll}
       1,  & \pi\msk1 j/q < x < \pi(j+1)/q,  \\
               
        0, & \mathrm{otherwise}, \\
     \end{array}  
   \right. 
\end{equation}
is the characteristic function of that subinterval.
In order to find the explicit values of the $a_{j}$, we first compute the Fourier coefficients of $\sigma^{j,q}(x)$:
\Eq{cjqk}
$$c^{j,q}_{k}=\mcases{
    \dsty   \frac{\i  \left(e^{-\i \pi k/q}-1\right)}{2\pii k }e^{-\i\pi j k/q},  & k\neq 0,  \\
      \dsty  \frac{1}{2\:q}, & k=0. \\
     } $$
Thus, by linearity, the Fourier coefficients of~\eq{upwc} are
\Eq{ajq}
$$c_k=\sum_{j=0}^{2q-1} \> a_{j}c^{j,q}_{k}.$$
Let us define the rescaled Fourier coefficients
\Eq{chatksum}
$$  \chat_k=\frac{1}{2\:q}\sum_{j=0}^{2q-1} \>a_{j}e^{-\i \pi \:j\: k/q} = 
     \mcases{ \frac{1}{2},  & k=0,  \\
               0,  & 0\neq k~\equiv 0\ \ \modx 2\:q,    \\
 \dsty              \frac{\pi \:k \:c_k}{i q \left(e^{-\i \pi\frac{\ell}{q}}-1\right)}, & k\not\equiv 0\ \ \modx 2\:q .}  $$
In particular, $\chat_0=1/2$, which implies
$$\sum_{j=0}^{2\:q-1} \>a_{j}=q.$$
For the problem at hand, we have
\begin{equation}  
   \chat_k = \mcases{   \frac{1}{2},  & k=0, \\
\dsty               \frac{e^{-\i \pi k^3p/q}}{ q \left(e^{-\i \pi\ell/q}-1\right)}, & {k\not\equiv 0}\ \ \modx 2\:q, \quad  k \roh{ odd}. } 
\end{equation} 
Inverting~\eq{chatksum} yields
\Eq{ajeq}
$$a_{j}=\sum_{\ell=-q}^{q-1}\;\chat_\ell\, e^{\fracz{\i \pi\ell j}{q}}=\frac{1}{2}+\frac{1}{q}\>\Re{\sum_{\substack{\ell=-q \\ \ell\ \textrm{odd}}}^{q-1}\;
\frac{e^{\fracz{-\i \pi\ell^3p}{q}}\,e^{\fracz{\i \pi\ell j}{q}}}{e^{\fracz{-\i \pi\ell}{q}}-1}},$$
{and hence}
\Eq{aj1}
$$a_{j}-a_{j-1}=-\,\frac{1}{q}\>\Re{\sum_{\substack{\ell=-q \\ \ell\ \textrm{odd}}}^{q-1}\;\zeta_{2q}^{j\ell-p\ell^3}   },$$
where, for conciseness, we adopt the notation
\Eq{zeta}
$$\zeta _m = e^{\fracz{2\pii \i}{m}}$$
for the primitive $m^\textrm{th}$ root of unity.

Now we define $\elt =\ell+q$ and consider the following cases:

\textbf{Case 1:} For $q$ odd,~\eq{aj1} becomes
$$a_{j}-a_{j-1}=\frac{(-1)^{p+j-1}}{q}\Re\left(\sum_{\substack{\elt =0 \\ \elt \ \textrm{even}}}^{2q-1}\zeta_{2q}^{j\:\elt  -p\:\elt \:^3}\right)=
\frac{(-1)^{p+j-1}}{q}\Re\left(\sum_{\nu=0}^{q-1}\;\zeta_{q}^{j\nu-4p\nu^3}\right),$$
which we rewrite as
\Eq{ajS1}
$$\boxed{a_{j}-a_{j-1}=\frac{(-1)^{p+j-1}}{q} \Se(p,q,j),}$$
where
\Eq{S1}
$$\Se(p,q,j)=\sum_{\nu=0}^{q-1}\zeta_{q}^{j\nu-4p\nu^3}.$$
The sum~\eq{S1} is real. Indeed, taking the complex conjugate yields
$$\overline{\Se(p,q,j)}=\sum_{\nu=0}^{q-1}\;\zeta_{q}^{-j\nu+4p\nu^3}=
\sum_{\kappa=0}^{q-1}\;\zeta_{q}^{\kappa j-4p\kappa^3}=\Se(p,q,j)$$
as the sum is over the complete residue class system $\ro{mod}\  q$.

\begin{remark}
Note that when $t=0$, \ro(and therefore $p=0$\ro), 
$$a_{j}-a_{j-1}=(-1)^{j-1}\frac{1}{q} \sum_{\nu=0}^{q-1}\;\zeta_{q}^{j\nu}=0, \soq{for} j \ne 0, q.$$
This result reconfirms that the initial condition only has jumps at $j=0$ and $j=q$. 
\end{remark}

\begin{remark}
We observe that, for $q$ odd, 
\begin{equation}
a_{j+q}-a_{j+q-1}=(-1)^{q}(a_{j}-a_{j-1})=-(a_{j}-a_{j-1}),
\end{equation}
which reconfirms our earlier observation that the jumps are symmetrically paired with those of opposite magnitude  around $j=q$. In addition, 
\begin{equation}
a_{j+q}+a_{j}=a_{j-1+q}+a_{j-1}, \, \soq{for all}j,
\end{equation}
which implies that $a_{j+q}+a_{j}=c$, with $c$ independent of $j$. Invoking~\eq{ajq}, we see that $c=1$.
\end{remark}

\textbf{Case 2:} For $q$ even,~\eq{aj1} becomes
$$
\eeq{
a_{j}-a_{j-1}=\frac{(-1)^{p+j-1}}{q}\Re\left(\sum_{\substack{\elt =0 \\ \elt \ \textrm{odd}}}^{2q-1}\zeta_{2q}^{j\:\elt  -p\:\elt \:^3}\right)=
\frac{(-1)^{p+j-1}}{q}\Re\left(\sum_{\substack{\kappa=0 \\ \kappa\ \textrm{even}}}^{2q-2}\zeta_{2q}^{j(\kappa+1)-p(\kappa+1)^3}\right)\\=
\frac{(-1)^{p+j-1}}{q}\Re\left(\zeta_{2q}^{j-p}\sum_{\nu=0}^{q-1}\;\zeta_{2q}^{2j\nu-2p\nu(4\nu^2+6\nu+3)}\right)\\=
\frac{(-1)^{p+j-1}}{q}\Re\left(\zeta_{2q}^{j-p}\sum_{\nu=0}^{q-1}\;\zeta_{q}^{j\nu-p\nu(4\nu^2+6\nu+3)}\right),}$$
which we rewrite as
\Eq{ajS2}
$$\boxed{a_{j}-a_{j-1}=\frac{(-1)^{j-1+p}}{q}\So(p,q,j),}. $$
where
\Eq{S2}
$$\So(p,q,j)=\zeta_{2q}^{j-p}\sum_{\nu=0}^{q-1}\zeta_{q}^{\nu (j-3p)-6p\nu^2-4p\nu^3} $$
The fact that~\eq{S2} is real will be demonstrated below. 

\begin{remark}
We observe that, for $q$ even, 
\begin{equation}
a_{j+q}-a_{j+q-1}=(-1)^{q+1}(a_{j}-a_{j-1})=-(a_{j}-a_{j-1})
\end{equation}
again confirming the fact that the jumps are symmetrically paired. 
\end{remark}

{Consider the following \emph{Kummer sum},~\rf{HBPa}:}
\Eq{TS}
$$
S(p,q,j)=\sum_{\nu=0}^{2q-1} \zeta_{2q}^{j\nu-p\nu^3} \roq{with} p\nmid q.
$$
One verifies, by taking the complex conjugate, that $S(p,q,j)\in\R$. Moreover, 
\begin{equation}
S(p,q,j)=\Se(p,q,j)+\So(p,q,j)=\sum_{\substack{\nu=0 \\ \nu\ \textrm{even}}}^{2q-2} \zeta_{2q}^{j\nu-p\nu^3}+\sum_{\substack{\nu=1 \\ \nu\ \textrm{odd}}}^{2q-1} \zeta_{2q}^{j\nu-p\nu^3}.
\end{equation}
Since $\Se(p,q,j)\in\R$, we conclude that $\So(p,q,j)\in\R$.  We will refer to $\Se(p,q,j)$ and $\So(p,q,j)$ as the \emph{even and odd partial Kummer sums}.  

As a final remark, consider the delta function initial condition $u(0,x)=\delta (x)$ that produces the fundamental solution $u(t,x)$ for the problem at hand. At the time value $t=\pi\msk1p/q$, the solution is
\begin{equation}
u(\pi\msk1 p/q,x)=\frac{1}{2\pii}\sum_{\nu=-\infty}^{\infty}e^{\i( \nu x-\pi  \nu^3p/q)}.
\end{equation}
Recalling the formula
\begin{equation}
\delta(x-\pi\msk1\ell/q)=\frac{1}{2\pii}\sum_{\nu=-\infty}^{\infty}e^{-\pi \i \ell \nu  /q}e^{ \i  \nu  x},
\end{equation}
we see that the solution is a finite linear combination of delta functions, namely
\begin{equation}
u(\pi\msk1 p/q,x)= \sum_{\ell=0}^{2q-1} {\beta_{\ell}(p/q)}\,\delta(x-\pi\msk1\ell/q) = 
\frac{1}{2\pii}\>\sum_{\ell=0}^{2q-1}\sum_{\nu=-\infty}^{\infty}{\beta_{\ell}(p/q)}\,e^{-\pi \i \ell \nu  /q}e^{ i\  \nu  x},
\end{equation}
provided
\begin{equation}
{\beta_{\ell}(p/q)}=\sum_{\nu=0}^{2q-1}\;e^{\pi \i( \nu \ell/q -  \nu^3  p/q)}=\sum_{\nu=0}^{2q-1}\;\zeta_{2q}^{ \nu \ell - p \nu^3 }.
\end{equation}
Thus, the coefficient ${\beta_{\ell}(p/q)}$ is given by the Kummer sum~\eq{TS}: 
\Eq{betal}
$${\beta_{\ell}(p/q)} = S(p,q,\ell).$$
{This is an important result in that the Kummer sums appearing in the fundamental solution contain all the information about partial Kummer sums which arise in the study of points of constancy.} Of course,
this relation is expected as once one knows the fundamental solution, one can obtain any other solution, and {hence number-theoretic Kummer sums play a fundamental} role in the underlying structure of \emph{any} solution to the periodic initial-boundary value problem for the Airy equation.

\section{Interlude: permutation polynomials}\label{sec:squarefree}

In this section, we recall the important concept of a permutation polynomial, as well as introduce some notation needed for the subsequent analysis. 
A \emph{permutation polynomial} over a 
 \emph{finite} ring is a polynomial which acts as a permutation of the elements of the ring. Here we only concern ourselves with the finite ring $\Z_q$ and polynomials of degree 3. 

A procedure to determine whether a polynomial has the permutation property can be found in~\rf{ChRyTa}. Consider the cubic polynomial
\Eq{p3}  
$$f(x)=f_1 x+f_2 x^2+f_3 x^3$$
over the ring $\Z_q$. If the conditions laid out in Table~\ref{tab:prime} for each prime $\qi$ appearing in the 
prime factorization~\eq{qfac} of $q$ and its corresponding exponent $n_i$ are satisfied, then we conclude that~\eq{p3} is a permutation polynomial, meaning that the elements of the ring $\Z_q$ are permuted among themselves by $f$.
\begin{table}
\begin{center}
\begin{tabular}{ |p{3cm}||p{2cm}|p{7cm}|  }
 \hline
 \multicolumn{3}{|c|}{Coefficient Test modulo $\textfrak{q}^n$} \\
 \hline
 $\textfrak{q}=2$ & $n=1$ &$(f_1+f_2+f_3)$ is odd. \\
     &$n>1$& $f_1$ is odd, $f_2$ is even, and $f_3$ is even.\\
 \hline
 $\textfrak{q}=3$ & $n=1$ &$(f_1+f_3)\neq0,\ f_2=0\,\ \modx 3$  \\
     &$n>1$& $f_1\neq0, \ f_1+f_3\neq0, \ f_2=0\,\ \modx 3$.\\
 \hline
 $3\mid \textfrak{q}-1$ & $n=1$ &$f_1\neq0,\ f_2=f_3=0\,\ \modx \textfrak{q}$  \\
             & $n>1 $  & $f_1\neq0, \ f_2=f_3=0\,\ \modx \textfrak{q}$  \\
 \hline
 {$3\nmid \textfrak{q}-1$} & $n=1$ &$ f_2^2=3f_1f_3, \ f_3\neq0\,\ \modx \textfrak{q} $  \\
  &  &\quad {\rm or} \ \ $f_1\neq0, \ f_2=f_3=0\,\ \modx \textfrak{q}$  \\
              &$n>1$   & $f_1\neq0, \ f_2=f_3=0\,\ \modx \textfrak{q}$  \\
  \hline
\end{tabular}
\end{center}
\caption{A test for conditions on primes $\qi$ for a permutation polynomial. Note that the second line for $3\nmid \textfrak{q}-1$ and $n=1$ will be irrelevant for our purposes.  
}\label{tab:prime}
\end{table}

\begin{remark}
For our purposes, when $q$ is odd and square free, the polynomial of interest will be $f(\nu)=-4p\nu^3+j\nu$, and therefore $f_1=j$, $f_2=0$ and $f_3=-4p$. 
The values are all $\ro{mod}\ q$, depending on the time value under consideration; see Section~\ref{sec:qodd}. 
When $q$ is even and square free, the polynomial of interest will be $g(\nu)=-4p\nu^3-6p\nu^2+(j-3p)\nu$ and therefore $g_1=j-3p$, $g_2=-6p$ and $g_3=-4p$. 
By demanding that $f$ or $g$ be a permutation polynomial, we obtain conditions on $j$, which in turn will correspond to points of constancy. 
\end{remark}

\begin{remark}
{Since permutation polynomials are only useful in our investigation when $q$ is square free, that is only the case $n=1$ in Table~\ref{tab:prime}, we omit the exponent $n$ from here on.}
\end{remark}
From now on, we denote the set of permutation polynomials at prime $\q$ as $\mathfrak{I}(q)$.
Note that if  $q$ has factorization \eq{qfac} with prime factors $\mathfrak{q_{i}}$, then the corresponding set of cubic permutation polynomials decomposes accordingly:
\begin{equation}
\mathfrak{I}(q)=\bigcap_{i}\;\mathfrak{I}(\mathfrak{q_i}).
\end{equation}

\section{The case when $q$ is odd }\label{sec:qodd}

In this section we assume $q$ is odd, and hence so are all its prime factors $\qi$.
Let us define 
\Eq{Pf}
$$P(x) = P(x\mid q,f_3)=\sum_{\nu=0}^{q-1} \;x^{f(\nu)}, \where f(\nu)=f_1 \nu +f_3\nu^3.$$
Since $f(\nu)$ is  defined $\ro{mod}\ q$, we have $\deg P(x) \leq \max f(\nu) \leq q-1$. In what follows, $f_1=j$ and, at times, we record the cubic coefficient $f_3$ in the definition of the 
polynomial $P(x)= P(x\mid q,f_3)$. 
Evaluating $P(x\mid q,f_3)$ at $x=\zeta_{q}$, we deduce the following identification 
\Eq{PS1}
$$P(\zeta_{q}\mid q,f_3)=\Se(p,q,j) .$$
Therefore, in view of~\eq{ajS1}, if $j$ is a point of constancy, $\zeta_{q}$ must be a root of $P(x)$.

Now we apply the Chinese Remainder Theorem in the following form reminiscent of the the Prime-Factor Algorithm (also known as Good-Thomas algorithm)~\rf{Good, TH}. Let us assume 
$q=N_1N_2$ with $\gcd(N_1,N_2)=1$. Consider the bijective re-indexing of the summation of~\eq{Pf} as 
\begin{equation}
\nu \longmapsto a\:N_1+b\:N_2, \where a\in[0,N_2-1], \quad b\in[0,N_1-1].
\end{equation}
Then,
\begin{equation}
f(\nu)\equiv f(a\:N_1)+f(b\:N_2) \mod (q = N_1N_2).
\end{equation}
Substituting this into~\eq{Pf}, we obtain
\begin{equation}
\ieq{ P(x\mid q,f_3)=\sum_{\nu=0}^{q-1} x^{f(\nu)} = \sum_{a=0}^{N_2-1} x^{f(aN_1)} \sum_{b=0}^{N_1-1} x^{f(bN_2)} \\
 =\sum_{a=0}^{N_2-1} (x^{N_1})^{f_1a+ f_3N_1^2a^3}\sum_{b=0}^{N_2-1} (x^{N_2})^{f_1 b+ f_3N_2^2b^3}
 =P(x^{N_1}\mid N_2,f_3N_1^2)\>P(x^{N_2}\mid N_1,f_3N_2^2),}
\end{equation}
noting that the cubic coefficient of $f(\nu)$ has changed for each factor, which is why we record it in our notation for $P(x)$. More generally, consider the prime factorization~\eq{qfac} of $q$. We 
let $\qhi$ denote the \emph{complement} of the prime factor $\qi$, which is the product of all factors \emph{except} for $\qi^{n_i}$. Repeating the previous process, one can show 
\Eq{Pfac}
$$P(x\mid q,f_3)=\prod_{i=1}^{m}\ P_{i}(x^{\qhi}\mid \qi^{n_i},f_3  \,\qhi^2).$$
Now, evaluating at $x=\zeta_{q}$, we have
\Eq{48}
$$P(\zeta_{q}\mid q,f_3)=\prod_{i=1}^{m}\ P_{i}(\zeta_{\qi^{n_i}}\mid \qi^{n_i},f_3  \,\qhi^2).$$
From~\eq{ajS1} and~\eq{PS1} this is equivalent to 
\Eq{49}
$$\boxed{\PC\left(\frac{p}{q}\right)=\bigcup_{i}\; \PC\left(\frac{p}{\qi^{n_i}}\right),}$$
which is the main structural result of this paper. It allows us to separately study the points of constancy at each prime $\qi$ in the prime factorization of $q$. In the following subsections we will 
consider various cases of $q$.  
\begin{remark}
Note that each prime $\qi$ is associated with its own polynomial $f(\nu)$ henceforth denoted as $f_{\qi}(\nu)$.
\end{remark}

\subsection{$q$ an odd prime}

We deal first with the easiest case, when $q$ is an odd prime. Let 
\begin{equation}
\Phi_{q}(x)=x^{q-1}+x^{q-2}+...+x+1
\end{equation}
denote the $q$-th cyclotomic polynomial, with $\Phi_{q}(\zeta_q)=0$. Since we are assume $q$ is an odd prime, $\Phi_{q}(x)$ is irreducible over $\mathbb{Q}$, and is the minimal polynomial of the primitive root of unity $\zeta_{q}$.

\Pr3
Assume $q$ is an odd prime. Let $f(\nu)=j\nu-4p\nu^3$.  If $ P(\zeta_{q})=0$, then $f(\nu)$ is a permutation polynomial $\ro{mod}\ q$. 
Conversely, if $f(\nu)$ is a permutation polynomial $\ro{mod}\ q$, then $\Phi_{q}(x)\mid P(x)$ and $j$ is a point of constancy. We denote this fact as
\begin{equation}
\Phi_{q}(x)\mid P(x) \quad\Longleftrightarrow\quad f(\nu)\in \mathfrak{I}(q).
\end{equation}

\begin{proof}
$(\Longrightarrow)$. If $\zeta_{q}$ is a root, this means that its minimal polynomial
$\Phi_{q}(x)$ will divide $P(x)$. Now, $\deg \Phi_{q}(x)=q-1$ and $\deg P(x)\leq q-1$. This means that $\Phi_{q}(x)=P(x)$ since $\Phi_{q}(1)=P(1)=q$. So, $f(\nu)$ has to be a permutation polynomial.
$(\Longleftarrow)$. If the polynomial $f(\nu)=j\nu-4p\nu^3$ assumes all values $\ro{mod}\  q$, then $P(\zeta_{q})$ is a 
sum of all the $q^\textrm{th}$ roots of unity, and hence is zero. As a result, $\Phi_{q}(x)\mid P(x)$ and we have a point of constancy at $j$.  
\end{proof}

Therefore, demanding that the polynomial $f(\nu)=j\nu-4p\nu^3$ be a permutation polynomial  for a certain prime $\q$ reduces Table~\ref{tab:prime} to Table~\ref{tab:reducedprime}.
\begin{table}
\begin{center}
\begin{tabular}{ |p{3cm}||p{2cm}|p{7cm}|  }
 \hline
 \multicolumn{3}{|c|}{Conditions for 
 $\PC$'s for all $p$ for prime $\q$} \\
 \hline
  $\textfrak{q}=3$ & $n=1$ &$j\not\equiv p \mod 3. $ \\
 \hline
$3\mid \textfrak{q}-1$ & $n=1$ &$\emptyset$  \\
\hline
$3\nmid \textfrak{q}-1$ & $n=1$ &$ j\equiv 0 \mod \q $  \\
  \hline
\end{tabular}
\end{center}
\caption{A reduction of Table~\ref{tab:prime} in the case $f(\nu)$ is a permutation polynomial.}\label{tab:reducedprime}
\end{table}
{For example}, at the particular time values $t=\pi/3,\ 2\pii/3,\ \pi/5,\ \pi/7,\ \pi/11$, and $\pi/13$ we find
$$\eeq{\PC\left(\frac{1}{3}\right)=\left\{j\in[0,5]\cap\Z\mid j\equiv 0,2\mod 3\right\} = \{0,2,3,5\},\\
\PC\left(\frac{2}{3}\right)=\left\{j\in[0,5]\cap\Z\mid j\equiv 0,1\mod 3\right\} = \{0,1,3,4\},\\
\PC\left(\frac{1}{5}\right)=\left\{j\in[0,9]\cap\Z\mid j\equiv 0\mod 5\right\}= \{0,5\},&
\PC\left(\frac{1}{7}\right)=\emptyset,\\
\PC\left(\frac{1}{11}\right)=\left\{j\in[0,21]\cap\Z\mid j\equiv 0\mod 11\right\}= \{0,11\}.}$$

\subsection{$q$ odd, composite, square-free}

When $q$ is square-free, every prime factor occurs just once in its prime factorization
\begin{equation}
q=\prod_{i=1}^{m}\qi,
\end{equation}
so $n_i=1$ for all $i$ in~\eqas{48}{49}.
It is important to note that the condition to be a permutation polynomial  for a specific prime $\qi$ involves the cubic polynomial 
\Eq{fqi3}
$${f_{\qi}(\nu)=j \nu -4\:p\:\qhi^2\nu^3.}$$
Here, expression \eq{49} allows us to consider each prime separately, each one giving rise to their own congruences. We then take the union of the individual congruences. 

To see this in practice, let us 
apply the conditions of the Table~\ref{tab:prime} for each selected prime $\qi$. First, suppose $\qi=3$. We easily see that $j\not\equiv p\,\qhi^{\:2}\mod3$. 
In the case $3\mid\qi-1$, we observe the condition $f_{3}\equiv 0 \mod \qi$ implies $-4p\,\qhi^{\:2}\equiv 0 \mod \qi$ but this is impossible. Therefore,
these primes do not contribute any congruences. Finally, for $3\nmid\qi-1$ we have 
\begin{equation}
0\equiv 3j(-4p)\qhi^{\:2}\ \ \modx \qi \quad\Longrightarrow\quad j\equiv 0 \ \ \modx \qi.
\end{equation}
All in all,
\Eq{oddqc}
$$\boxed{\PC\left(\frac{p}{q}\right)=\left\{j\in[0,2q-1]\cap\Z\mid  (j\not\equiv p\,\qhi^2\ \modx3) 
 \cup (j\equiv 0\ \modx \qi \mid 3\nmid\qi-1) \right\}.}$$

Here are a few examples: 
$$\eeq{\PC\left(\frac{1}{15}\right)=\left\{j\in[0,29]\cap\Z\mid (j\equiv 0,2\ \modx 3)\cup (j\equiv 0\ \modx 5) \right\},\\
\PC\left(\frac{2}{15}\right)=\left\{j\in[0,29]\cap\Z\mid (j\equiv 0,1\ \modx 3)\cup (j\equiv 0\ \modx 5) \right\},\\
\PC\left(\frac{1}{21}\right)=\left\{j\in[0,41]\cap\Z\mid j\equiv 0,2\ \modx 3 \right\}.}$$
   
\subsection{$q$ an odd prime power}

Assume that $q$ is an odd number with prime factorization~\eq{qfac}. As above, given a prime factor $\qi^n$, let $\qhi$ be its complement in $q$.
First, we consider the case $n=2$, namely\footnote{{The $\Vert$ notation means that $\mathfrak{q}_{i}^2$ \emph{fully divides} $q$, that is, $\mathfrak{q}_{i}^2 \mid q $  but $\mathfrak{q}_{i}^3 \nmid q $.}} $\mathfrak{q}_{i}^2 \Vert q $. For the case $n=3$, we will explicitly consider the time value $t=\pi p/3^3$ as we shall encounter congruences $\mod 3^2$. For higher 
powers, the treatment is conceptually similar although computationally much more involved.

{As in \eq{fqi3},} consider the polynomial $f_{\qi}(\nu)=j\nu-4\:p\:\qhi^2\nu^3\mod q$, where from now on we suppress the index $i$. 
Consider
the following parametrization for $\nu$ as $\nu=a+b\:\q$ with $a,b\in\mathbb{F}_{\q}$ and therefore $\nu\in\mathbb{F}_{\q}\times\mathbb{F}_{\q}$. Now $f_{\q}(\nu)$ 
becomes
\begin{equation}
f_{\q}(\nu)=f_{\q}(a+b\:\q)=aj-4\:p\:\qh^2a^3+b(j-12a^2p\:\qh^2)\q.
\end{equation}
Note that as $f_{\q}(\nu)$ is defined mod $q$, both summands are as well, therefore $(j-12\:a^2p\:\qh^2)$ makes sense mod $\q$. Now, the relevant factor of $P(x)$ in~\eq{Pfac} is
\Eq{Psq}
$$\eeq{P_{\q}(x)=\sum_{\nu=0}^{\q^2-1}\; x^{f_{\q}(\nu)}=\sum_{a=0}^{\q-1}\sum_{b=0}^{\q-1}\; x^{ja-4p\:\qhat^2a^3+b(j-12\:a^2p\:\qhat^2)\q}\\ =
\sum_{a=0}^{\q-1}x^{ja-4p\qhat^2a^3}\sum_{b=0}^{\q-1} \;(x^\q)^{b(j-12\:a^2p\:\qhat^2)}=\sum_{a=0}^{\q-1}\;x^{f(a)}\sum_{b=0}^{\q-1} (x^\q)^{b(j-12a^2p\:\qhat^2)}.}
$$

In order to see how $\q$ contributes to the factorization 
\eq{Pfac}, we ask whether $x=\zeta_{\q^2}$ is a root of the polynomial~\eq{Psq}:  $P_{\q}(\zeta_{\q^2})=0$. We have
\begin{equation}
P_{\q}(\zeta_{\q^2})=\sum_{a=0}^{\q-1}\;\zeta_{\q^2}^{f_{\q}(a)} \ \sum_{b=0}^{\q-1} \;(\zeta_{\q^2}^\q)^{b(j-12\:a^2p\:\qhat^2)}=
\sum_{a=0}^{\q-1}\;\zeta_{\q^2}^{f(a)}\ \sum_{b=0}^{\q-1}\; \zeta_{\q}^{b(j-12\:p\:\qhat^2a^2)}.
\end{equation}
For this to vanish, we need $j-12\:p\:\qhat^2a^2\not\equiv 0 \mod\q$  for all $a\in \mathbb{F}_{\q} $. If $\q=3$, then we immediately obtain $j\not\equiv 0 \mod 3$. For other primes, we wish to 
state this more elegantly, using the  \emph{Legendre symbol}~\rf{HaWr}, namely
\begin{equation}
\left(\frac{12^{-1}p^{-1}(\qhat^2)^{-1}j}{\q}\right)_2=-1 \quad\Longrightarrow\quad \left(\frac{3\:p\:j}{\q}\right)_2=-1.
\end{equation}
Finally, the contribution of a prime $\q$ for which $\q^2\Vert q$ to $\PC\left({p}/{q}\right)$ is 
\Eq{oddq2}
$$
   \boxed{\PC\left(\frac{p}{\q^2}\right)= \mcases{
       j\not\equiv 0 \ \ \modx 3,  && \q=3,  \\
    \ostrut{17}{17}           \left(\frac{j}{\q}\right)_{2}=1, & & \left(\frac{3\:p}{\q}\right)_{2}=-1,    \\
               \left(\frac{j}{\q}\right)_{2}=-1, && \left(\frac{3p}{\q}\right)_{2}=1. }}
$$

Consider the example $p=1, \q=7^2=49,\ \q=7$. Clearly,
\begin{equation}
\left(\frac{3}{7}\right)_{2}=-1 \quad\Longrightarrow\quad  \left(\frac{j}{7}\right)_{2}=1 \quad\Longleftrightarrow\quad  j=1,2,4 \ \ \modx 7.
\end{equation}
As an additional example, when $p=1, q=5^2=25, \q=5$, it easily follows
\begin{equation}
\left(\frac{3}{5}\right)_{2}=-1 \quad\Longrightarrow\quad  \left(\frac{j}{5}\right)_{2}=1 \quad\Longleftrightarrow\quad  j=1,4  \ \ \modx 5.
\end{equation}

To proceed further, for $n>2$ and $\mathfrak{q}>3$, we can invoke the \emph{Stationary Phase Method} laid out in~\rf{Coch1,Coch2,Fish} and the references therein. The basic idea behind this method  is that exponential sums for prime powers can be computed by adapting the well known Method of Stationary Phase that applies to integrals. As when computing integrals, which requires  determining the relevant \emph{critical points}, the analogous notion of \emph{critical point congruence} carries over to the case of exponential sums. 
\emph{It is important to note that outside the critical point congruence the exponential sum is zero\/}! This property is precisely described by \eq{oddq2}, being the content of Theorem 1.1 in~\rf{Coch2}.

As a final consideration, we will investigate the case $n=3$ with $\mathfrak{q}=3$, whereby $\qh=1$, corresponding to the time value $t=\pi \:p/3^3=\pi \:p/27$. As before, we write 
$\nu=a+b\:\mathfrak{q}+c\:\mathfrak{q}^2 $ with $a,b,c \in\mathbb{F}_{\mathfrak{q}}$ and therefore $\nu\in\mathbb{F}_{\mathfrak{q}}\times\mathbb{F}_{\mathfrak{q}}\times\mathbb{F}_{\mathfrak{q}}$. Hence
$f(\nu)=\nu j-4\:p\:\nu^3\ \modx q$ becomes 
\begin{equation}
f(\nu)=f(a+b\:\mathfrak{q}+c\:\mathfrak{q}^2)=aj-4\:p\:a^3+b(j-12\:a^2p)\mathfrak{q}+(c\:j-12\:p\:a(a\:c+b^2))\mathfrak{q}^2.
\end{equation}
Now, $P(x)$ becomes 
\begin{equation}
P(x)=\sum_{a=0}^{\mathfrak{q}-1}\;x^{f(a)} \ \sum_{b=0}^{\mathfrak{q}-1} \;(x^\mathfrak{q})^{b(j-12\:a^2p)}\ \sum_{c=0}^{\mathfrak{q}-1}\; (x^{\mathfrak{q}^2})^{c\:j-12\:p\:a(a\:c+b^2)},
\end{equation}
and, evaluating at $x=\zeta_{\mathfrak{q}^3}$, we obtain 
\begin{equation}
P(\zeta_{\mathfrak{q}^3})=\sum_{a=0}^{\mathfrak{q}-1}\;\zeta_{\mathfrak{q}^3}^{f(a)}\ \sum_{b=0}^{\mathfrak{q}-1} \;\zeta_{\mathfrak{q}^2}^{b(j-12\:p\:a^2)} \ 
\sum_{c=0}^{\mathfrak{q}-1}\; \zeta_{\mathfrak{q}}^{c\:j-12\:a\:p(a\:c+b^2)},
\end{equation}
which could also be written as 
\begin{equation}
P(\zeta_{\mathfrak{q}^3})=\sum_{a=0}^{\mathfrak{q}-1}\;\zeta_{\mathfrak{q}^3}^{f(a)}\ \sum_{b=0}^{\mathfrak{q}-1}\; \zeta_{\mathfrak{q}^2}^{b(j-12\:p\:a^2)}\zeta_{\mathfrak{q}}^{-12\:p\:a\:b^2}\ 
\sum_{c=0}^{\mathfrak{q}-1}\; \zeta_{\mathfrak{q}}^{c(j-12\:p\:a^2)}.
\end{equation}

For the case $\mathfrak{q}=3$, we have $j-12\:p\:a^2\not\equiv0\mod 3$ for all $a$ and therefore $j\equiv1,2 \mod 3$. This was the congruence $\mod 3$. In order to obtain the congruence $\mod 3^2$,
we look at the sum over $b$ and we observe that it simplifies because $\zeta_{3}^{-12\:p\:a\:b^2}=1$ for all $b\in [0,2]$. We 
demand $j-12\:p\:a^2\equiv0\mod 3$ for all $a$ and therefore $j= 3\:k$. The sum over $b$ becomes 
\begin{equation}
\sum_{b=0}^{2} \zeta_{3}^{b(k-4\:p\:a^2)}
\end{equation}
and we require $k-4\:p\:a^2\not\equiv0\mod 3$ for all $a$ which, applying the Principle of Quadratic Reciprocity for the prime $\mathfrak{q}=3$, yields $k\equiv2\:p\mod 3$ and therefore $j\equiv6\:p\mod 3^2$. All in 
all,
\Eq{oddq3}
$$\boxed{\PC\left(\frac{p}{3^3}\right)=\left\{j\in[0,53]\cap\Z\mid (j\equiv 1,2 \ \modx3)\cup (j\equiv 6\:p \ \modx3^2) \right\}.}
$$

\section{The case when $q$ is even} \label{sec:qeven}

When $q$ is even, we use both versions of the odd partial Kummer sum~\eq{S2} implied by~\eq{cjqk}:
\begin{equation}
\So(p,q,j)=\sum_{\substack{\elt =0 \\ \elt \ \textrm{odd}}}^{2q-1}\zeta_{2q}^{j\:\elt  -p\:\elt \:^3}=\zeta_{2q}^{j-p}\sum_{\nu=0}^{q-1}\;\zeta_{q}^{j\nu-p\nu(4\nu^2+6\nu+3)}.
\end{equation}
When $q$ is a power of $2$, we employ the first sum, while when $q$ is even composite, we use the second version. Note that $p$ must be an odd number. 

\subsection{$q$ a power of $2$} 

Consider the case $q=2^n$.
We define the polynomial
\begin{equation}
Q(x)=\sum_{\substack{\nu=1 \\ \nu\ \textrm{odd}}}^{2q-1} x^{f(\nu)}.
\end{equation}
Clearly,  when $f(\nu)=j\nu-p\:\nu^3 \ \modx 2q$, we have
\begin{equation}
Q(\zeta_{2q})=\sum_{\substack{\nu=1 \\ \nu\ \textrm{odd}}}^{2q-1} \zeta_{2q}^{j\nu-p\nu^3}=\So(p,q,j).
\end{equation}
Since $f(\nu)\leq 2q-1$, this implies that $\deg Q(x)\leq 2q-1$. Therefore, if $j$ is a point of constancy, $\zeta_{2q}$ must be a root of $Q(x)$. 

Consider first the simplest case with $n=1$, so $q = 2$, corresponding to the time value $t=\pi \:p/2$ and therefore $j\in[0,3]$. 
We have $\nu=a+2b$ with $a\in\mathbb{F}_{2}^{\times}$, $b\in\mathbb{F}_{2}$ and therefore $\nu\in\mathbb{F}_{2}^{\times}\times\mathbb{F}_{2}$. Thus 
\begin{equation}
f(\nu)=f(a+2b)=aj-pa^3+2b(j-3a^2p)=j-p+2b(j-3p)\mod 4,
\end{equation}
as $a=1$. Hence
\begin{equation}
Q(x)=\sum_{\substack{\nu=1 \\ \nu\ \textrm{odd}}}^{3} x^{f(\nu)}=\sum_{b=0}^{1} \;x^{j-p+2b(j-3p)}=x^{j-p}\sum_{b=0}^{1} \;(x^{2})^{b(j-3p)}.
\end{equation}
We now evaluate 
\begin{equation}
Q(\zeta_{4})=\zeta_{4}^{j-p}\sum_{b=0}^{1} (\zeta_{4}^{2})^{b(j-3p)}=\zeta_{4}^{j-p}\left(1+(-1)^{j-3p}\right) =0    \quad\Longleftrightarrow\quad j\equiv 0\mod2.
\end{equation}
We conclude that
\Eq{even2}
$$\boxed{\PC\left(\frac{p}{2}\right)=\left\{j\in[0,3]\cap\Z\mid j\equiv 0\ \modx 2  \right\} = \{0,2\}.}
$$

We consider the case with $n=2$, so $q = 2^2 = 4$, corresponding to the time value $t=\pi \msk1p/4$ and therefore $j\in[0,7]$. 
We have $\nu=a+2b+4c$  with $a\in\mathbb{F}_{2}^{\times}$ and $b,c\in\mathbb{F}_{2}$, hence $\nu\in\mathbb{F}_{2}^{\times}\times\mathbb{F}_{2}\times\mathbb{F}_{2}$. Thus
\begin{equation}
f(\nu)=f(a+2b+4c)=j-p+2b(j-3p)+4(cj-3pc-3pb^2)\mod 8,
\end{equation}
as we have already set $a=1$. Thus,
$$\eeq{{Q(x)= \sum_{\substack{\nu=1 \\ \nu\ \textrm{odd}}}^{7} x^{f(\nu)}}=\sum_{b=0}^{1} \sum_{c=0}^{1}\;x^{j-p+2b(j-3p)+4(cj-3pc-3pb^2)} \\
=x^{j-p}\sum_{b=0}^{1} \;(x^{2})^{b(j-3p)}(x^{4})^{-3pb^2}\sum_{c=0}^{1}\; (x^{4})^{c(j-3p)}.}
$$
Evaluating at $x=\zeta_{2q}=\zeta_{8}$ yields
\begin{equation}
Q(\zeta_{8})=\zeta_{8}^{j-p}\sum_{b=0}^{1}\; (\zeta_{8}^{2})^{b(j-3p)}(\zeta_{8}^{4})^{-3pb^2}\sum_{c=0}^{1}\; (\zeta_{8}^{4})^{c(j-3p)}=\zeta_{8}^{j-p}\left(1+(-1)^p\zeta_{4}^{j-3p}\right)\left(1+(-1)^{j-3p}\right).
\end{equation}
As $p$ is  odd, we have 
\Eq{even4}
$$\boxed{\PC\left(\frac{p}{4}\right)=\left\{j\in[0,7]\cap\Z\mid (j\equiv 0 \ \modx 2)\cup (j\equiv 3p \ \modx4) \right\}.}
$$

As {our final case, suppose} $n=3$, so $q = 2^3 = 8$,  corresponding to the time value $t=\pi \msk1p/8$ and therefore $j\in[0,15]$. We 
have $\nu=a+2b+4c+8d$  with $a\in\mathbb{F}_{2}^{\times}$ and $b,c,d\in\mathbb{F}_{2}$. Thus,\begin{equation}
f(\nu)=f(a+2b+4c+8d)=j-p+2b(j-3p)+4(cj-3pc-3pb^2)+8(dj-pb^3-pd)\mod 16.
\end{equation}
Hence
\begin{equation}
 \begin{split}
Q(x)= & \sum_{\substack{\nu=1 \\ \nu\ \textrm{odd}}}^{15} x^{f(\nu)}=\sum_{b=0}^{1} \sum_{c=0}^{1}\sum_{d=0}^{1}\;x^{j-p+2b(j-3p)+4(cj-3pc-3pb^2)+8(dj-pb^3-pd)} \\
=&x^{j-p}\sum_{b=0}^{1}\; (x^{2})^{b(j-3p)}(x^{4})^{-3pb^2}(x^{8})^{-pb^3}\ \sum_{c=0}^{1}\; (x^{4})^{c(j-3p)}\ \sum_{d=0}^{1} \;(x^{8})^{d(j-p)}.
\end{split}
\end{equation}
We now evaluate at $x=\zeta_{2q}=\zeta_{16}$ to obtain
\begin{equation}
\begin{split}
Q(\zeta_{16})=& 
\zeta_{16}^{j-p}\left(1+(-1)^p\zeta_{4}^{-3p}\zeta_{8}^{j-3p}\right)\left(1+\zeta_{4}^{j-3p}\right)\left(1+(-1)^{j-p}\right).
\end{split}
\end{equation}
As $p$ is necessarily odd, we obtain the congruences $j\equiv0\mod2$, $j\equiv (2-p)\mod4$ and $j\equiv (4-2p)\mod8$. The last congruence is subsumed in the first one. All in all,
\Eq{even8}
$$\boxed{\PC\left(\frac{p}{8}\right)=\left\{j\in[0,15]\cap\Z\mid (j\equiv 0 \ \modx 2)\cup (j\equiv (2-p) \ \modx 4) \right\}.}
$$

\begin{remark}
One might be tempted to think that the factorization of $Q(\zeta_{8})$ or $Q(\zeta_{16})$ persists in higher powers of $2$. This is not the case. The crucial difference between $2$ and the odd primes is that in the former case we always set the constant term $a=1$. 
\end{remark}

For higher powers of $2$, experimental evidence with \textsc{Mathematica} demonstrates that
\begin{itemize}
  
\item $q = 16$, \ $t=\pi/16$: 
In this case, we require that $j$ be even or {$j\not\equiv7 \mod 8$}.
\item $q = 32$, \ $t=\pi/32$: 
In this case, we require that $j$ be even or {$j\not\equiv7 \mod 8$}.
\end{itemize}

\subsection{$q$ even, composite}

Let $q$ be even and composite, so $\qi = 2$ occurs as a factor.
Consider the polynomials 
\begin{equation}
Q(x\mid q,g_2,g_3)=\sum_{\nu=0}^{q-1}\; x^{g(\nu)},
\end{equation}
where $g(\nu)=g_1 \nu +g_2\nu^2+g_3\nu^3$ is defined $\ro{mod}\  q$,
hence $\deg Q(x)\leq q-1$. In what follows, $g_1=j-3p$, $g_2=-6p$, and $g_3=-4p$.   

Given this notation, evaluating $Q(x\mid q,f_3)$ at $x=\zeta_{q}$ we have the following identification with the odd partial Kummer sum:
\begin{equation}
Q(\zeta_{q}\mid q,g_2,g_3)=\So(p,q,j) .
\end{equation}
Therefore, if $j$ is a point of constancy, $\zeta_{q}$ must be a root of $Q(x)$. 
As before, we assume $q=N_1N_2$. Consider the {summation reparametrization} 
\begin{equation}
\nu \longmapsto aN_1+bN_2, \qquad a\in[0,N_2-1], \quad b\in[0,N_1-1],
\end{equation}
so that
\begin{equation}
g(\nu)\equiv g(aN_1)+g(bN_2) \mod  (q = N_1N_2).
\end{equation}
Substituting back into~\eq{Pf}, we obtain
\begin{equation}
 \begin{split}
 Q(x\mid q,g_2,g_3)&=\sum_{\nu=0}^{q-1} x^{g(\nu)} = \sum_{a=0}^{N_2-1} x^{g(aN_1)} \sum_{b=0}^{N_1-1} x^{g(bN_2)} \\
 &=\sum_{a=0}^{N_2-1} (x^{N_1})^{g_1a+g_2a^2N_1+g_3N_1^2a^3}\sum_{b=0}^{N_2-1} (x^{N_2})^{g_1a+g_2a^2N_2+g_3N_2^2a^3}\\
 &=Q(x^{N_1}\mid N_2,g_2N_1,g_3N_1^2)\>Q(x^{N_2}\mid N_1,g_2N_2,g_3N_2^2).
  \end{split}
\end{equation}
Thus, in general,
\begin{equation}
Q(x\mid q,g_2,g_3)=\prod_{i=1}^{m} \;Q_{i}(x^{\qhi}\mid \qi^{n_i},g_2 \qhi, g_3 \qhi^2).
\end{equation}
Note that both the cubic coefficient \emph{as well as the quadratic} coefficient of $g(\nu)$ change for each prime factor.
Now, evaluating at $x=\zeta_{q}$
\begin{equation}
Q(\zeta_{q}\mid q,g_2,g_3)=\prod_{i=1}^{m}\; Q_{i}(\zeta_{\qi^{n_i}}\mid \qi^{n_i},g_2 \qhi, g_3 \qhi^2),
\end{equation}
which leads us to the decomposition
\Eq{49e}
$$\boxed{\PC\left(\frac{p}{q}\right)=\bigcup_{i} \;\PC\left(\frac{p}{\qi^{n_i}}\right)},$$
with one of the $\qi$'s being $2$. This generalizes condition~\eq{49} to include even numbers.

\subsection{$q$ even, composite, square free}

For the case where $q$ is even and square free, with prime factorization~\eq{qfac}, 
one of the prime factors is $\qi = 2$, and $p$ must be odd. We consider now the family of polynomials 
\begin{equation}
g_{\qi}(\nu)=(j-3p) \nu -6p\:\qhi\nu^2 -4p\:\qhi^2\nu^3, 
\end{equation}
 indexed by a selected prime $\qi$.

We apply the conditions of the Table~\ref{tab:prime} (notice that we relabel $f_{i}\mapsto g_{i} $) for each selected prime $\qi$. So, for $\qi=2$, we demand that 
$j-3p+(-6p)\qhi+(-4p)\qhi^2$ be odd. This immediately implies $j$ must be even for all $p$. Now consider $\qi=3$. We easily obtain $j\not\equiv p\:\qhi^2\mod3$. 
Moving on to the case when $3\mid\qi-1$, we observe that the condition $g_{2}\equiv 0 \mod \qi$ implies $-6\:p\:\qhi\equiv 0 \mod \qi$ but this is impossible. Finally,
for the case $3\nmid\qi-1$ we have 
\begin{equation}
(-6p)^2\qhi^2\equiv 3(j-3p)(-4p)\qhi^2 \mod \qi \quad\Longrightarrow\quad j\equiv 0 \mod \qi.
\end{equation}
All in all,
\Eq{evenqn}
$$\boxed{\ibeq{62}{\PC\left(\frac{p}{q}\right)=\left\{\ostrut{12}{7}j\in[0,2q-1]\cap\Z\ \right| \\\left.\ostrut{12}{7}(j\equiv 0 \ \modx 2)\cup (j\not\equiv p\:\qhi^2\ \modx3)
{} \cup (j\equiv 0 \ \modx \qi \mid 3\nmid\qi-1) \right\}.}}
$$
 Let us consider some examples. When $q = 6$ and $10$, corresponding to  the time values $t=\pi/6$ and $\pi/10$, we  obtain
\begin{equation}
{\eeq{\PC\left(\frac{1}{6}\right)=\left\{j\in[0,11]\cap\Z\mid (j\equiv 0 \ \modx 2)\cup (j\equiv 0,2 \ \modx3 ) \right\},\\\PC\left(\frac{1}{10}\right)=\left\{j\in[0,19]\cap\Z\mid (j\equiv 0 \ \modx 2)\cup (j\equiv 0 \ \modx5 ) \right\}.}}
\end{equation}

\subsection{The general even case}

Assume that $q$ is an even number and its prime factorization is given in~\eq{qfac}. Let us choose a certain prime $\qi^n$ in that factorization and let $\qhi$ be its complement in $q$.
{Here we will} deal only with the case $n=2$. We consider the polynomial 
$$f_{\qi}(\nu)= (j-3\:p)\nu-6\:p\:\qhi\nu^2-4\:p\:\qhi^2\nu^3\ \ \modx q,$$ 
and from now on suppress the index $i$. 

We carry out $\textit{mutatis mutandis}$ the analysis of the previous section and we end up with
\begin{equation}
f_{\q}(\nu)=f_{\q}(a+b\:\q)=a(j-3\:p)-6\:p\:\qhat a^2-4\:p\:\qhat^2a^3+b((j-3\:p)-12\:a\:p\:\qhat-12\:a\:p^2\qhat^2)\q,
\end{equation}
and
\begin{equation}
\eeq{P_{\q}(x)=\sum_{a=0}^{\q-1\;}x^{a(j-3\:p)-6\:p\:\qhat a^2-4p\qhat^2a^3}\ \sum_{b=0}^{\q-1}\; (x^\q)^{b((j-3\:p)-12\:a\:p\qhat-12\:a\:p^2\qhat^2)}
\\=\sum_{a=0}^{\q-1}\;x^{f_{\q}(a)}\ \sum_{b=0}^{\q-1} \;(x^\q)^{b((j-3\:p)-12\:a\:p\qhat-12\:a\:p^2\qhat^2)}.}
\end{equation}
Evaluating at $x=\zeta_{\q^2}$, we obtain
\begin{equation}
P_{\q}(\zeta_{\q^2})=\sum_{a=0}^{\q-1}\;\zeta_{\q^2}^{f_{\q}(a)}\ \sum_{b=0}^{\q-1}\; 
(\zeta_{\q^2}^\q)^{b((j-3\:p)-12\:a\:p\:\qhat-12\:a\:p^2\qhat^2)}=
\sum_{a=0}^{\q-1}\;\zeta_{\q^2}^{f_{\q}(a)}\ \sum_{b=0}^{\q-1}\; \zeta_{\q}^{b((j-3\:p)-12\:a\:p\qhat-12\:a\:p^2\qhat^2)}.
\end{equation}
We need
\begin{equation}\label{eq:leg}
(j-3\:p)-12\:a\:p\qhat-12\:a\:p^2\qhat^2\not\equiv 0 \mod\q \quad\Longrightarrow\quad  j-3\:p(1+2\:a\:\qhat)^2\not\equiv 0 \mod\q,
\end{equation}
 for all $a\in \mathbb{F}_{\q} $. If $\q=3$, then we immediately obtain $j\not\equiv 0 \mod 3$. For other primes, multiplication by any number coprime to $\q$ is an
automorphism of $\mathbb{F}_{\q}$ fixing $0$ and translating by $1$ is again an automorphism of $\mathbb{F}_{\q}$. For $\q\neq3$, (~\ref{eq:leg}) can be stated using the Legendre symbol as
\begin{equation}
 \genfrac(){}{0}{3\:p\:j}{\q}_{2}=-1 \quad\Longrightarrow\quad \genfrac(){}{0}{3\:p}{\q}_{2} \genfrac(){}{0}{j}{\q}_{2}=-1.
\end{equation}
Finally, the contribution of a prime 
$\q$ for which $\q^2\Vert q$, to $\PC\left(\dsty\frac{p}{q}\right)$ is 
\Eq{evenq2}
$$   \boxed{\PC\left(\frac{p}{\q^2}\right)= \left\{
     \begin{array}{lrl}
       j\not\equiv 0\ \ \modx 3,  && \q=3,  \\\dsty
              \ostrut{22}{22} \left(\frac{j}{\q}\right)_{2}=1,  && \dsty\left(\frac{3\:p}{\q}\right)_{2}=-1,     \\
              \dsty \left(\frac{j}{\q}\right)_{2}=-1, && \dsty\left(\frac{3\:p}{\q}\right)_{2}=1. \\
     \end{array}   
   \right. }
$$

\section{Conclusion}\label{sec:conclusion}

We have considered the one-dimensional periodic linearized Korteweg--deVries equation, also known as the Airy partial differential equation, subject to two different initial conditions: a delta function situated at the origin, leading to the 
fundamental solution, and a piecewise constant step profile. Upon temporal evolution, at rational times relative to the overall length of the interval, in the former case the fundamental solution takes 
the form  of a \emph{Dirac comb}. In the latter case, the solution remains piecewise constant on suitable subintervals. The coefficients of the individual delta functions in the Dirac comb profile as well as 
the difference between two successive constant values in the piecewise constant profile are both given by number-theoretic Kummer sums.

The goal of this work was to understand the points of constancy  of the piecewise constant profiles and the ``absence'' of delta functions in the Dirac comb profiles. Both of these translate into the requirement 
that a particular  Kummer sum vanish. While only the Kummer sums corresponding to the piecewise constant profile have been analyzed in detail,  our methods work equally well for any Kummer sum. We find that, at a given rational time  value $t=\pi \:p/q$, each prime power $\q^n$ in the prime factorization of the denominator $q$ gives rise to  certain congruences $\modx \q^{i}$ for $1\leq i < n$, which, in 
addition, depend on the numerator $p$. The easiest prime to deal with is $\q=2$, where we are able to go up to $n=3$. In cases when $\q$ is an odd prime and $n=2$, we made use of the  $\textit{Principle of Quadratic Reciprocity}$ to perform the required analysis. As noted above, higher powers require a more in depth analysis as the required  computations become increasingly intricate. 
However, when $\q=3$, we have been able to resolve the power $n=3$.  Higher powers and other primes can potentially be analyzed using the Stationary Phase Method for  exponential sums described in~\rf{Coch1,Coch2,Fish}.

We would also like to mention some questions which merit further investigation. First, it would be important to analyze the relationship between the Kummer sums associated with the piecewise
constant profiles with their ``$x$-derivatives'', meaning their respective delta function configuration. Second, one may consider the following question:
Where are the maximal or minimal jumps in the piecewise constant profile or, equivalently, the maximal and minimal contributions to the Dirac comb? 

One ultimate goal of these studies is to understand the phenomenon of ``fractalization" at irrational times using insights from such number theoretic analysis. For the linearized Korteweg--deVries equation, fractalization was first established in~\rf{Osk1}; recent results and bounds on fractal dimensions can be found in \rf{ChErTz,ErSh}. More 
specifically, one needs to consider sequences of rational time values approaching an irrational time and how the discontinuous but piecewise constant profiles converge to the fractal, non-differentiable but continuous profile. Such sought-after results may well depend on how closely one can approximate the irrational time value by rationals, that is,  on its \emph{irrationality measure}, \rf{Roth}.  This quest will involve dealing with high prime powers rendering it quite challenging. Nevertheless, we hope that the approach taken in this work may eventually shed some light on this problem.  

{Further, as noted in \rf{BerryKlein,BMS} in the case of the periodic linear Schr\"odinger equation, the entire solution graph $u(t,x)$ forms an intricate non-differentiable surface exhibiting anisotropic behavior --- meaning differing fractal profiles when restricted to the lines in the $t\:x$ plane, e.g., for varying $t$ and fixed $x$, or lines having other slopes --- whose detailed understanding is far from complete. 
Recent deep results of Erdo\u gan and Shakan, \rf{ErSh} now provide bounds on these various fractal dimensions.  On the other hand, determining the corresponding more informative Hausdorff dimensions seems beyond current analytical capabilities.} 

\boh{Acknowledgements:}
The first author would like to thank Adrian Diaconu for advice and references.  The second author would like to thank Professor Vasilis Papageorgiou and Professor Aristides Kontogiorgis for enlightening discussions and especially Professor Pavlos Tzermias for enlightening discussions as  well as inspiration 
during the course of many years. Both authors also wish to thank Natalie Sheils for her remarks and suggestions.  {We also thank the referees for their careful reading of the paper, corrections, and suggestions.}

\boh{Data accessibility:} The paper contains no experimental or computational data other than what is listed in the tables therein which can be easily reproduced using the formulas and methods described.

\boh{Competing interests:}  We have no competing interests.

\boh{Authors' contributions:} Both authors contributed to the results and the writing of the paper, and both gave final approval for publication.

\boh{Funding statement:} The work was not supported by external funding.

\boh{Ethics statement:} The work is purely mathematical and involves no experiments.

\end{document}